\documentclass[12pt]{article}

\usepackage{amssymb,url, amsthm, amsmath, verbatim}

\newtheorem{theorem}{\noindent Theorem}

\newtheorem{definition}{\noindent Definition}
\newtheorem{corollary}{\noindent Corollary}

\newtheorem{remark}{\noindent Remark}

\author{F.~Petrov, A.~Vershik}

\date{30.06.09}

\title{Uncountable Graphs and Invariant Measures on the Set
of Universal Countable Graphs}

\begin{document}
\maketitle

\begin{abstract}
We give new examples and describe the complete lists of all measures on the set of
countable homogeneous universal graphs and $K_s$-free homogeneous universal graphs
(for $s\geq 3$) that are invariant with respect to the group of all permutations
of the vertices. Such measures can be regarded as random graphs
(respectively, random $K_s$-free graphs). The well-known example of
Erd\"os--R\'enyi (ER) of the random graph corresponds to the Bernoulli measure on the set of
adjacency matrices. For the case of the universal $K_s$-free
graphs there were no previously known examples of the invariant measures on the space of such graphs.

The main idea of our construction is based on the new notions of
{\it measurable universal}, and {\it topologically universal} graphs, which are interesting themselves.
The realization of the construction can be regarded as two-step randomization for universal measurable graph
: {\it "randomization in vertices"} and {\it "randomization in edges"}. For $K_s$-free, $s\geq 3$
there is only randomization in vertices of the measurable graphs. The completeness of our lists
is proved using the important theorem by D.~Aldous about $S_{\infty}$-invariant matrices,
which we reformulate in appropriate way.

\end{abstract}

\tableofcontents

\section{Introduction: problem and results}

\subsection{Universal graphs}

Fix a countable set $V$ and consider the set ${\cal G}_V$ of all
graphs (undirected, without loops and multiple edges) with
$V$ as the set of vertices. Equip ${\cal G}_V$ with the weak
topology (the base of the weak topology is formed by the collections of sets of graphs
that have a given induced graph structure on a given finite set of vertices).
The weak topology allows us to define the notion of
Borel sets and Borel $\sigma$-field on ${\cal G}_V$, and to consider Borel probability measures on
$\cal{G_V}$. It is convenient to take the set of positive integers $\Bbb N$
as $V$.

We can identify a graph $\Gamma$ with its adjacency matrix
$A_{\Gamma}$:  an entry $e_{i,j}$, $i,j \in \Bbb N$, of
$A_{\Gamma}$ is equal to $1$ or $0$ if $(i,j)$ is an edge or not an
edge, respectively. Thus the space ${\cal G}_{\Bbb N}$ of graphs can be identified with
the space $M^{Sym}_{\Bbb N}(0;1)$ of all infinite symmetric zero-one matrices with
zeros on the principal diagonal, equipped with the usual weak topology
on the space of matrices.

The infinite symmetric group ${\frak S}^{\Bbb N}$ of all
permutations of
the set $\Bbb N$ acts naturally on the space of graphs ${\cal{G}}_{\Bbb N}$.
Each orbit of ${\frak S}^{\Bbb N}$ is a class of
isomorphic graphs, and the stabilizer  of a given graph, as a subgroup
of ${\frak S}^{\Bbb N}$, is the group of all automorphisms of this graph.
The action of ${\frak S}^{\Bbb N}$ is continuous with respect to the weak topology
on ${\cal{G}}_{\Bbb N}$, and to
the weak topology on the group ${\frak S}^{\Bbb N}$ itself. In terms of
the space of matrices $M^{Sym}_{\Bbb N}(0;1)$, this action obviously means a
simultaneous permutation of the rows and columns of the adjacency
matrices. The action naturally extends to an action on Borel
measures on the spaces of graphs and matrices.

We will consider subsets of ${\cal G}_{\Bbb N}$ that
are invariant under the action of ${\frak S}^{\Bbb N}$, and {\it
invariant Borel probability measures} on such sets. Of most interest
are subsets of ${\cal G}_{\Bbb N}$ on which the group ${\frak S}^{\Bbb N}$
acts transitively; namely, an example important for our purposes is
the
family of universal graphs in a category of graphs.

Consider a small category $\cal C$ whose objects are finite or
countable graphs (the sets of vertices of these graphs are subsets of
${\Bbb N}$) that contains a {\it universal object}. This means that
there is an object of $\cal C$, a countable graph $\Gamma$, that satisfies
the following properties:

1) $\Gamma$ contains any finite graph of the
category $\cal C$  as a subobject (up to isomorphism)

and

2) the group of all isomorphisms of $\Gamma$
acts transitively on the set of isomorphic finite subgraphs of
$\Gamma$.

Such graphs are called homogeneous universal graphs of the
category  $\cal C$. Hereafter we just call them ``universal'', without
explicit mentioning homogeneity.
Fra{\"\i}ss\'e's theory (see, e.g., \cite{W}) gives
transparent necessary and sufficient conditions for the
description of categories that have a universal graph. We may
assume that the sets of vertices of all universal graphs in all these
categories coincide with the whole set ${\Bbb N}$, so we can identify graphs with
their adjacency matrices from $M^{Sym}_{\Bbb N}(0,1)$, and the set of universal
graphs is an orbit of the action of the group ${\frak S}^{\Bbb N}$. {\it By
a ``random graph'' in a given category we mean a ${\frak
S}^{\Bbb N}$-invariant Borel probability measure on the set of graphs
that is concentrated on the set of universal graphs of this category.}

Here we restrict ourselves to the following category: ${\cal C}_s$,
$s>2$, is the category of all {\it finite or countable graphs that
contain no $s$-cliques $K_s$} (an $s$-clique is a complete graph with $s$ vertices,
$s>2$). Also denote by $\cal C$ the category of all finite or
countable graphs. It is well known that Fra{\"\i}ss\'e's axioms are valid in these
cases, and there are universal graphs in all these categories.
A corollary of the existence of universal graphs asserts that all
universal graphs are mutually isomorphic, so a universal graph is
unique up to isomorphism; consequently, the set of all universal graphs
is an orbit of the group ${\frak S}^{\Bbb N}$.
We describe the set of invariant measures on these orbits.

\subsection{Random graphs and invariant measures on the set of the universal graphs}
 We consider a ``random $K_s$-free graph,'' which is the same as a
 ${\frak S}^{\Bbb N}$-invariant measure on the set of universal $K_s$-free graphs.
The existence of a ${\frak S}^{\Bbb N}$-invariant measure on the set of ordinary
 universal graphs (the category $\cal C$) is well known: this is the
 Erd\"os--R\'enyi \cite{ER} random graph. In our terms, the examples
 of Erd\"os and R\'enyi are
 the Bernoulli measures on the space $M^{Sym}_{\Bbb N}(0,1)$ of adjacency
 matrices with the distribution $(p,1-p)$, $0<p<1$, for each entry.
 Note that for $p=1/2$ this Bernoulli measure is the weak limit of
 the uniform  measures on the sets of finite graphs with $n$
 vertices as $n$ tends to infinity. We will see that there are many
 other ${\frak S}^{\Bbb N}$-invariant measures on the set of universal graphs.

 As to the categories ${\cal C}_s$, $s>2$, no invariant measures (or no
 random graphs) were known at all. For the case $s=3$, it was known
 that the weak limit of the uniform measures on the set of finite
 triangle-free graphs with $n$ vertices as $n$ tends to infinity is
 not a measure on the set of universal triangle-free graphs, but an invariant
 measure on the set of universal bipartite graphs. This follows from
 the results of \cite{EKR,KRP} on asymptotic estimations of the
 number of odd cycles in typical triangle-free graphs\footnote{We are grateful to
 Professor G.~Cherlin for the references
 to these papers.}. This means that the uniform measure, as an approximation
 tool, is too rough for obtaining the desired measure. {\it Nevertheless we
 proved that there exist uncountably many invariant ergodic measures on the set of
 $K_s$-free graphs for $s>2$.}
  Note the paradoxical fact that, in spite of the transitivity of the action of
 the group ${\frak S}^V$ on the set of universal graphs, there
 exist uncountably many different (pairwise singular)
 ${\frak S}^{\Bbb N}$-invariant ergodic measures; this is a new manifestation
 of Kolmogorov's effect, see details  in \cite{V1}.

 Remark that our construction of the universal continuous graph for the case $s=2,3$ is {\it
shift invariant} which means that there is a transitive action of a group $\Bbb R$ on the set of vertices (which is $\Bbb R$) of the
continuous graph. For the countable universal graph the existence of the
transitive action of the group $\Bbb Z$ on the set of vertices
is trivial; for the case of triangle free universal graph it was proved by C.~Henson \cite{He}; who also had proved nonexistence of
such action for $s>3$. \footnote{We are grateful to the reviewer who pointed out to this paper.}
We also mentioned that fact for continuous case.

\subsection{How uncountable universal graphs can help toward countable ones:
double randomization.}

 For constructing ${\frak S}^V$-invariant measures on the space of
universal graphs, we will use a very natural general method of constructing
invariant measures on the set of infinite matrices.
It looks like the Monte-Carlo or randomization method.
Specifically, we take a {\it continuous graph}, that is, a standard measure space $(X,m)$,
regarded as the set of vertices, and a subset $E\subset X^2$, regarded as the set of edges,
and then choose vertices (points of $X$) at random,
independently, with distribution $m$; the induced countable subgraph is our random graph.
If we want to obtain an invariant measure on the set of
universal ($K_s$-free universal, etc.) graphs, we must use (and first define!)
a universal (respectively, $K_s$-free universal, etc.) continuous graph.
{\it Thus our examples of invariant measures on the space of universal graphs
come from ``randomization in vertices'' of universal continuous measurable
graphs.}\footnote{Note that our notion of universality of continuous graphs
is not a categorical universality and homogeneity.}
Note that the notion of a universal continuous graph is perhaps
of interest in itself in the theory of models and ``continuous combinatorics.''
It looks similar to the universal Urysohn space if we compare it with the countable
universal metric space. We will consider this analogy in a separate
paper.\footnote{The notion of a continuous graph in general must be very useful in variational
 calculus, geometrical optimal control, etc.}

The method of the randomization of the vertices does not give all invariant measures on the set of universal graphs
(or on the set of all countable graphs). Even Erdos-Renji example of random graph is not of that type.
In order to describe all invariant measures on the set of universal graphs as well as invariant measures on the other sets of the countable graphs, we must generalize this method and use {\it another kind of randomization, namely, ``randomization in edges''}.
In this paper we shortly describe this construction based on the notion of {\it generalized graph} and
on the important theorem due to D.~Aldous \cite{A}, which describes in some sense all ${\frak S}^{\Bbb N}$-invariant measures
on the space of the infinite matrices. In particulary we apply
this construction for the universal graphs. It gives the list of
all ${\frak S}^{\Bbb N}$-invariant ergodic measures on the set
of universal or $K_s$-universal graphs.

Remark that in order to prove that our constructions
exhaust the list of all invariant measures on the set of
universal ($K_s$-free universal) graphs we use the important theorem due
to D.~Aldous \cite{A} about invariant measures on the
space of matrices. We formulate that theorem in a suitable
version, which will be considered with a new proof of it
by the second author in the separate place. See also \cite{V}), where these problems are linked to
the problem of classification of measurable functions of
the several variables.

Thus our scheme looks like the following transitions:
 universal Borel graph with measures $\rightarrow$
 topologically universal graph ($\rightarrow$ homogeneous
 topologically universal
 graph)\footnote{The homogeneity is used for the ordinary and triangle-free cases only.} $\rightarrow$
 randomization in vertices $\rightarrow$ invariant measures on the set of countable universal graphs
 $\rightarrow$ randomization in edges $\rightarrow$ the list of all invariant measures on the set of countable universal
 ($K_s$-free universal) graphs. In brief, our description of
 invariant measures reduces to the choice of a deterministic continuous
 graph, then to randomization of its vertices (randomly choosing some vertices), and then to
 randomization of edges.

 The method of this paper does not help to solve the problem due to Prof. G.~Cherlin
 about existence of the finite triangle free "almost" universal graph. The reason is
 that it is difficult to extract from our constructions the implicit type of
 finite dimensional approximations of the constructed invariant measures.
 In the same time the continuous models for constructions of random countable
 objects can be applied in many situations.

\subsection{About this paper.}
 Let us give a short description of the contents of the paper.
 In the second section we consider the notions of continuous graphs
 and universal continuous graphs of various types using a
 generalization of the criterion of universality. We give two kinds of
 definitions: for measurable (Borel) graphs and topological graphs;
 the latter ones are more convenient for our goals.
 Section~3 is devoted to a particular construction of topologically universal
 ($K_s$-free universal) graphs. We define even
a shift-invariant graph structure with $\Bbb R$ as the set of vertices for the
 ordinary and triangle-free cases. This gives the existence of nontrivial
 ${\frak S}^{\Bbb N}$-invariant measures on the set of universal graphs.
 The main part of the Section~4 has deal with the general constructions of the
 invariant measures not only for universal graphs. We give the classification of measures
 obtained in terms of the randomization in edges in the spirit of paper \cite{V}.
 Then we define the generalized graph and give a general scheme of the double randomization
 of the universal continuous graphs (in vertices and edges). This is tightly connected with
 the mentioned above Aldous's theorem about the list of all ${\frak S}_{\Bbb N}$-invariant
 measures on the ${0-1}$ matrices.
  This gives a list of all such invariant measures for universal and $K_s$-free ($s>2$) universal graphs.
   Some problems and comments are collected in the last section.
 One of the main practical problems is to find directly the finite-dimensional
 distributions of our measures on the set of universal graphs, or, more specifically, to
 describe the approximation of random universal graphs in our sense in
 terms of random finite graphs.

  Professor T.~Tao informed the second author that the idea of
 using continuous graphs has already appeared in
 the recent papers by L.~Lovasz and his coauthors \cite{L1, L2}, where an
 analog of a continuous weighted graph was defined.
 In \cite{Diak}, this notion was also associated with Aldous' theorem.
 Our goals and constructions are different from those constructions:
 we consider {\it universal continuous graphs}.

The authors are grateful to Professors N.~Alon, G.~Cherlin and T.~Tao for important references,
the reviewer of the paper for  very useful comments
and to Prof. N.~Tsilevich for her help with preparing the final version of
the paper and useful criticism.

\section{Theme and variation on universal graphs}

\subsection{Countable graphs: the criterion of universality}

 Recall that the universality of a countable graph $\Gamma_u$ is equivalent
 to the following two conditions:

 $(i)$ any finite graph $\gamma$ can be isomorphically embedded into
 $\Gamma_u$;

 $(ii)$ for any two isomorphic finite induced subgraphs $\gamma_1$, $\gamma_2$
 of $\Gamma_u$, any isomorphism between them can be extended to an
 isomorphism of the whole graph $\Gamma_u$.

 It is easy to prove that the following well-known criterion is equivalent to $(i)\&(ii)$
 (see, e.g., \cite{Cam}):

 $(iii)$ for any two disjoint finite subgraphs $\gamma_1\subset \Gamma_u$
  (call it ``black'') and $\gamma_2 \subset \Gamma_u$ (call it ``white'') there exists a
 vertex $v\in \Gamma_u$ that is joined with the white vertices
 and is not joined with the black ones.

 Now we will give an analog of this condition for the case of
{\it  triangle-free and $K_s$-free graphs.}
\begin{theorem}
{\rm1.} A countable triangle-free graph $\Gamma$ is a universal triangle-free
graph if and only if the  following condition is satisfied:

$(iii_3)$ for any two disjoint finite subgraphs $\gamma_1\subset \Gamma_u$
  (call it ``black'') and $\gamma_2 \subset \Gamma_u$
 (call it ``white''), where the white subgraph has no edges, there exists a
 vertex $v\in \Gamma_u$ that is joined with all white vertices
   and is not joined with the black vertices.

{\rm2.} For $s>2$, a countable $K_s$-free graph $\Gamma$ is a universal $K_s$-free
graph if and only if the following condition is satisfied:

$(iii_m)$ For any two disjoint finite subgraphs $\gamma_1\subset
\Gamma_u$ (call it ``black'') and
$\gamma_2 \subset \Gamma_u$ (call it ``white''), where the white subgraph
is $K_{s-1}$-free,
there exists a
 vertex $v\in \Gamma_u$ that is joined with all white vertices
 and is not joined with the black vertices.
\end{theorem}

Of course, the first part of the theorem is a special case of the second one,
and in what follows we will consider the triangle-free case as a special case
of $K_s$-free graphs with $s=3$.
 The proof of the theorem is a simple modification
 of the proof of the
previous theorem.

\subsection{Universal measurable graphs}

Now we give the definition of Borel (measurable), topologically
universal, and topologically universal $K_s$-free graphs for $s>2$. But
first of all we will give the definition of continuous graphs themselves.
Our definitions of these notions are not of the greatest possible generality, but they
are appropriate for our goals.

 Recall that a standard (uncountable) Borel space $X$ is a space with
 a fixed $\sigma$-field of subsets that is Borel isomorphic to the
 interval $[0,1]$ equipped with the $\sigma$-field of Borel subsets.

\begin{definition}
A Borel (undirected)
graph is a pair $(X,E)$ where $X$ is a standard
Borel space and $E\subset X\times X$ is a symmetric Borel subset
in $X\times X$ that is disjoint from the diagonal $\{(x,x), x \in
X\}$.
\end{definition}

We will denote  $E_x=\{y\in X: (x,y)\in X\}$ and $E'_x=X\backslash
E_x$. Note that if $\{x_k\}_{k=1}^{\infty}$ is a finite or countable sequence
in $X$, then it can be regarded as an ordinary finite or countable
subgraph of $(X,B)$ with the induced graph structure.  We say that a Borel
graph is {\it pure} if $E_x\ne E_y$ for $x\ne y$. Note that universal countable
graphs are pure.

The following measure-theoretic definition is more  useful for us.
\begin{definition}
A measurable (Borel) graph is a triple
$(X,m,E)$ where $(X,m)$
is a standard Lebesgue space with a continuous probability measure $m$
(i.e., the pair $(X,m)$ is isomorphic in the measure-theoretic sense
to the interval $[0,1]$ equipped with the Lebesgue measure)
and $E \subset X\times X$ is a symmetric measurable
set of positive $(m \times m)$-measure.\footnote{Strictly speaking,
we must consider the class of sets
that are equal to $E$ up to a set of zero measure;
consequently, we define a class of ($\bmod\, 0$)-coinciding graphs.}
\end{definition}

A measurable graph is called {\it pure} if the measurable map $x \rightarrow E_x(\bmod 0)$
from $X$ to sigma algebra of $mod 0$-classes of measurable sets is injective $\bmod 0$.

\begin{definition}
A universal (respectively, $K_s$-free universal)
measurable graph is a pure measurable graph $(X,m,E)$ that satisfies the following
property:
for almost all sequences $\{x_k\}_{k=1}^{\infty} \in X^{\infty}$
with respect to the Bernoulli measure $m^{\infty}$ in the space
$X^{\infty}$,  the induced countable graph on the set of vertices
$\{x_k\}$ is a universal countable graph (respectively, a $K_s$-free
universal countable graph).\footnote{It is more correct to call
such graphs countably universal, because the condition deals
only with countable subsets of $X$.}
\end{definition}

The definition above is indirect, but it is not difficult to formulate
direct definition which is equivalent to the previous.
\begin{theorem}
{\rm1.}
The pure measurable graph $(X,m,E)$ is universal in the above sense iff
for almost any two disjoint
finite sets $\{x_1, \dots, x_n\}\in X$ and
$\{y_1, \dots, y_m\}\in X$ the $m$-measures of the
following intersections:
$$m(\bigcap_{i,j} (E_{x_i}\cap E'_{y_j}))$$
are positive;
{\rm2.} The pure measurable graph $(X,m,E)$ is $K_s$-free universal iff
there are almost no $s$-tuples in $X$ for which the induced (by the set
$E$) graph is a $K_s$-graph;
and for any positive integers $k,t$ and for almost any two finite subsets $x=\{x_1,\dots,
x_k\}$, $\{y_1, \dots, y_t\} \subset X$ such
that the induced graph $x$ has no $K_{s-1}$-subgraphs, the $m$-measure of the following
intersection is positive:
 $$m(\bigcap_{i,j} (E_{x_i}\cap E'_{y_j}))>0.$$

For $s=3$ this gives the definition of a triangle-free
topologically universal graph.
\end{theorem}

The proof of the equivalence consists in direct application of the ergodic theorem
(or even individual law of large numbers) to the indicators of intersections defined above.

A direct corollary of our definition is the following theorem
 which will be used in what follows.

\begin{theorem}[Construction of invariant measures]
Let $(X,m,E)$ be a universal
(respectively, $K_s$-free universal) measurable graph. Consider
the space
$$X^{\infty}=\prod_{n=1}^{\infty} (X,m)$$ and the map
$$F:X^{\infty}\rightarrow M_{\Bbb N}(0,1),$$
$$F(\{x_i\})=\{e_{i,j}\},\qquad e_{i,j}=\chi_E(x_i,x_j),\quad i,j \in \Bbb N,$$
where $\chi_E$ is the characteristic function of the set $E\subset
X\times X$. Denote by $F^*$ the map defined on the space of Borel measures
on $X^{\infty}$ by the following formula: if $\alpha$ is a Borel measure on  $X^{\infty}$, then
$[F^*(\alpha)](C)=\alpha(F^{(-1)}(C))$, $C \subset M_{\Bbb N}(0,1)$.
  Then the measure $F^*(m^{\infty})\equiv
\mu_{\{X,m,E\}}$ is an ${\frak S}_{\Bbb N}$-invariant measure on the set of universal (respectively,
$K_s$-free universal) homogeneous countable graphs.
\end{theorem}
\begin{proof}
Follows from the fact that the Bernoulli measure $m^{\infty}$
is ${\frak S}^{\Bbb N}$-invariant.
\end{proof}
The following formula gives the implicit expression of the measure of cylindric sets
(it does not use the condition of universality). Suppose  $A=\{a_{i,j}\}, i,j=1 \dots n$
is $(0-1)$-matrix of order $n$, and $C_A$ is a cylindric set of all infinite $(0-1)$-matrices, which has the matrix $A$ as submatrix
on the NW-corner. Then the value of the measure $\mu_{\{X,m,E\}}\equiv \mu$ on the cylinder $C_A$ is

$$\mu(C_A)=m^n\{(x_1,x_2 \dots x_n):(x_i,x_j)\in E; \quad \mbox{if}\quad  a_{i,j}=1; \quad (x_k,x_r) \notin E \quad  \mbox{if} \quad a_{k,r}=0\},$$
where $m^n=m\times \dots (n) \dots \times m$

Because the given criterion of the universality is difficult to check for measurable graphs,  we will
use topological version of universality which is much more convenient.

\subsection{Topologically universal graphs}

As we have mentioned, it is not easy to check that a given measurable graph $(X,m,E)$
is a universal measurable graph. For this reason,
we will give a more restrictive definition of {\it topological
universality}, whose conditions are easier to check.

Let us define a topologically universal graph. For simplicity, we assume
that $X$ is a Polish (= metric separable complete) space, but
this is not necessary.

Given a set $Y\subset X$, denote its complement by $Y'=X\setminus Y$ and its
closure by $\bar{Y}$. A {\it topological graph (undirected, without loops)
is a pair $(X,E)$ where $X$ is  a Polish space and $E\subset X\times
X$ is a closed subset that has a nonempty interior and an empty
intersection with the diagonal.}\footnote{Our definition allows
vertices to have uncountably many neighbors. There are
many other definitions of topological graphs and topological graphs
with weights; one of them uses the notion of a polymorphism or Markov
transformation.}
Put $E_x=\{y\in X: (x,y)\in X\}$.
We say that a topological graph is pure if $E_x\ne E_y$ for $x\ne y$.

\begin{definition}{\rm1.} A pure topological graph  $(X,E)$ is called
topologically universal if the set $E$ satisfies the following property:

$(U)$ For any two disjoint finite sets $\{x_1, \dots, x_n\}\in X$ and
$\{y_1, \dots, y_m\}\in X$, the intersection
$$\bigcap_{i,j} (E_{x_i}\cap E'_{y_j})$$ has a nonempty interior.

{\rm2.} A topological graph  $(X,E)$ is called topologically
$K_s$-free universal if

$(U_m)$ there are no $m$-tuples in $X$ for which the induced (by the set
$E$) graph is a $K_s$-graph (a complete $m$-graph);
and for any positive integers $k,t$ and any two finite subsets $x=\{x_1,\dots,
x_k\}$, $\{y_1, \dots, y_t\} \subset X$ such that the induced graph on $x$ has no
$K_{s-1}$-subgraphs, the set
 $$\bigcap_{i,j} (E_{x_i}\cap E'_{y_j})$$
has a nonempty interior in $X$.

For $s=3$ this gives the definition of a triangle-free
topologically universal graph.
\end{definition}

It is worth mentioning that a topologically universal graph is not a
``universal topological graph'' in the sense of the category of topological
graphs; our definition is more flexible. As in the case of measurable
graphs, it is more correct to call it a ``countably universal
topological graph.''

Recall that a Borel measure on a Polish space is called
nondegenerate if it is positive on all nonempty open sets.

\begin{theorem} Let $(X,E)$ be a topologically universal graph
(respectively, a topologically universal $K_s$-free graph, $s>2$); then for
every nondegenerate Borel probability measure $m$ on the space $X$,
the triple $(X,E,m)$ is a universal measurable (respectively,
universal measurable $K_s$-free) graph in the sense of the
definition of Section~{\rm2.2}.
\end{theorem}

\begin{proof}
Let $m$ be a nondegenerate measure on $X$. We must check that the
property $(iii)$ (respectively, $(iii_M)$) from Section~2.1 is
valid for almost all (with respect to the Bernoulli measure
$m^{\infty}$) sequences $\{x_k\}$. First of all, almost all
sequences $\{x_k\}$
are everywhere dense in the separable metric space  $X$.
Consequently, every such sequence   $\{x_k\}$ contains points from any
open set in $X$. Since all the sets $\cap_{i,j} E_{\dots}\cap E'_{\dots}$,
described in Definition 4, have a nonempty
interior, the proof is done. The measurable graph is pure since
the topologically universal graph is pure.
\end{proof}

Using this theorem, we immediately obtain the following corollary, which shows how to
produce required measures on the set of universal graphs.

\begin{corollary}
For every nondegenerate measure $m$ on a topologically universal
(respectively, triangle free, $K_s$-free) graph $(X,E)$,
 the  measure $\mu_{\{X,m,E\}}$  is a ${\frak S}^{\Bbb
N}$-invariant measure on the set of universal (respectively,
universal triangle-free, universal $K_s$-free) countable graphs.
\end{corollary}

 The existence of topologically universal graphs and topologically universal
$K_s$-free graphs is proved in the next section.

The reason why we introduce, in addition to the notion of a measurable
universal graph, the notion of a topologically universal graph
is that it is difficult to formulate a measure-theoretic analog
of the property that the interiors of the sets $E(x,y)$ are nonempty,
or equivalent properties, which are important for extending a
countable graph structure to a continuous one.
But there are no doubts that this notion is useful in itself.

\begin{remark}
{\rm All previous definitions can be written in a more rigid form if we use
the invariance with respect to an action of a group on the set
of vertices of the graph. Let $G$ be a group, and let
the set of vertices $X$ be a $G$-space; we can repeat our definitions of continuous
and universal graphs for a $G$-invariant graph structure. For example,
let $X=G$, and let the set of edges $E\subset G\times G$ be left
$G$-invariant: $E =\{(g,h): g^{-1}h \in Z\}$, where $ Z \subset G $.
Group-invariant graph structures (Cayley objects in the terminology of
 \cite{Cam}) were considered in \cite{Cam,CV,V2}.}
\end{remark}

\section{Construction of continuous homogeneous graphs}

 We will prove the existence of topologically universal graphs and
 topologically universal $K_s$-free graphs. According to the previous results, our construction gives examples
 of invariant measures on the space of universal graphs. As we will see, there are
 many such constructions which produce uncountably many nonequivalent
 invariant measures. We choose the simplest example, namely, consider
 the additive group $X=\mathbb{R}$ as the set of vertices of a topological
 graph and define an appropriate set of edges (a subset of $\mathbb{R}^2$).
 Moreover, for the case of ordinary universal graphs and triangle-free
 graphs, we suggest a graph structure that is {\it shift-invariant}:
 $$E=\{(x,y): |x-y|\in Z\}\subset\Bbb R^2,$$ where the set
 $Z\subset (0,+\infty)$ will be constructed by induction. This means
 that the additive group $\mathbb{R}$ acts by the automorphisms transitively on the set of
 vertices.

 We will prove the following main result.
\begin{theorem}
{\rm1.} There is a topologically universal graph (respectively,
topologically universal
triangle-free graph) with the additive group $\Bbb R$ as the
set of vertices and a shift-invariant graph structure.

{\rm2.} There is a topologically universal  $K_s$-free graph (for $s>3$) with the additive
group $\Bbb R$ as the set of vertices. There is no universal
$K_s$-free graph for $s>3$ with a shift-invariant graph structure.
\end{theorem}

\begin{proof}
1. We will begin with the construction of a countable universal
(triangle-free universal) graph with the additive group of
rational numbers  $\Bbb Q$ as the set of vertices and a shift-invariant graph
structure. After that we extend the construction onto  $\Bbb R$.

We choose $X=\mathbb{Q}$ as the set of vertices
 and construct a set  $Z\subset \mathbb{Q}$ that will be the subset of vertices
joined by edges with zero $\textbf{0} \in \mathbb{Q}$. Thus $(x,y)$
 is an edge of our graph if and only if $|x-y|\in Z$, $x\ne y$. We construct $Z$ as the union of
 disjoint nondegenerate intervals of $\mathbb{Q}$ with
 rational endpoints such that any bounded
 set $M\subset \mathbb{R}$ contains only finitely many such intervals.
  The required shift-invariant structure of a continuous universal
 (triangle-free universal) graph on  $\mathbb{R}$ will appear if $\bar Z$,
  the closure of $Z$, is the set of vertices $x \in \mathbb{R}$ that are joined by
  edges with $\textbf{0} \in \mathbb{Q}$. In a sense, it is a
 completion of that universal (triangle-free universal) rational graph.

 It is easy to reformulate the conditions of universality in terms
 of the set $Z \subset {\Bbb R}_+$ using the shift-invariance:

 For a universal graph, we obtain the following condition.

\smallskip
 $(U)$ For every pair of disjoint finite sets of rational
 numbers $\{x_1,\dots, x_k\}$, $\{y_1,\dots ,y_t\}$ there exists a rational
 number $c$ such that $|c-x_i| \in Z$, $i=1,\,2,\dots,k$; $|c-y_j|
 \notin Z$, $j=1,\,2,\dots,t$.
\smallskip

 For a universal triangle-free graph, the condition on the set $Z$
 is more rigid:

\smallskip
$(UTF)$ (a) The {\it sum-free condition}: The equation $x+y=z$ has no solutions with
$x, y, z \in \bar{Z}$ (this is a corollary of the triangle-free condition
for graphs).

 (b) for every pair of disjoint finite sets of rational numbers
 $\{x_1,\dots,x_k\}$, $\{y_1,\dots ,y_t\}$ such that  $|x_i-x_j|\notin
 {\bar Z}$, $1\le i<j\le k$, there exists a rational number $c$ such
 that $|c-x_i| \in Z$, $i=1,2,\dots,k$; $|c-y_j| \notin Z$,
 $j=1,2,\dots,t$.
\smallskip

 In both cases,
 our construction will satisfy a stronger condition, which is necessary
 for our purposes: there exists an interval $(c_1,c_2)$ of points
 $c$ satisfying the above property.

 The construction of the set $Z$ is inductive and based on the enumeration
 of arrays of points from  $\mathbb{Q}$. We will use the simplest
 method of enumeration suitable both for
 ordinary and triangle-free graphs.

 Choose $\gamma$, a  pair of finite sets of disjoint intervals;
 the first set $$\{(a_1,a'_1),(a_2,a'_2),\dots, (a_k,a'_k)\}=\gamma^a$$ of the pair
 will be called ``white,'' and the second set
 $$\{(b_1,b'_1),(b_2,b'_2),\dots, (b_\ell,b'_\ell)\}=\gamma^b$$ will be called
 ``black''; all the closures of these $k+\ell$ intervals are mutually disjoint.
 We will call such a pair $\gamma$ a pattern. In the triangle free case require that
 the closure of white part of the pattern is sum-free (i.e. consider
 only such patterns).
 There are countably many patterns,
 so we can label all patterns with positive integers
 $\gamma_1, \gamma_2, \dots $.
 Note that for every pair $\bar x=(x_1,\dots,x_p)$, $\bar y=(y_1,\dots, y_q)$
 of disjoint finite subsets of $\mathbb{Q}$,
  there exists a pattern whose white part
contains the set $\bar x$ and black part contains the set $\bar y$.

For each pattern $\gamma_n$ we will define by induction  a set $Z_n$, the union
of finitely many intervals with rational endpoints, such that for all $z\in Z_n$,  $u\in Z_{n+1}$
we have $z<u$ (the monotonicity condition). Then $Z$ will be the union of these $Z_n$: $Z=\cup_n Z_n$.

As the induction base we can take an arbitrary pattern, or even
empty set.

 Now we will consider two cases.

 1) Construction of a universal continuous graph.

 Assume that we have already constructed sets
 $Z_1,\dots,  Z_{n-1}$, each of which is the union
 of closed disjoint intervals and satisfies the
 monotonicity condition above. Assume also that the condition $(U)$ is satisfied
 for all patterns with numbers less than $n-1$. More exactly,  if a
 set $\bar x =(x_1, \dots, x_k)$ belongs to the white part of a pattern
 with number less than $n-1$ and a set $\bar y=(y_1,\dots, y_t)$
 belongs to the black part of this pattern,
 then there exists an open interval $C \subset \mathbb{Q}$ such that
 $|c-x_i|\in \cup_1^{n-1}
 Z_i$ and $|c-y_j|\notin \cup_1^{n-1} Z_i $  for
 every $c\in C$.
 Consider the next pattern
 $\gamma_n=\left(\gamma^a=\cup_1^k (a_i,a'_i), \gamma^b=\cup_1^s(b_j, b'_j)\right)$
 and define a set $Z_n$ as follows. Find such a large $c$
 that
 $$
 c > \max_{i,j}\{a_i, b_j\}+ \max\{z: z \in \cup_{i=1}^{n-1} Z_i\}+1,
 $$
 and put $Z_n=\cup_{i=1}^k(c-a'_i-\varepsilon,c-a_i+\varepsilon)$. It is clear that
 for small enough $\varepsilon$ all points
 that belong to a sufficiently small neighborhood of $c$ are joined by edges with the white part
 of $\gamma_n$ and are not joined with the black part of $\gamma_n$, because the shifted
 segments $[c-a'_i,c-a_i]$ and $[c-b'_j,c-b_j]$ $i=1,\dots, k$, $j=1,\dots, s$, are disjoint,
 and so their small neighborhoods are disjoint as well.
  This completes the construction of the set $Z=\cup_n Z_n$.
 Now let us check that the graph with the set of vertices $\mathbb{Q}$
 and the edges $\{(x,y):|x-y|\in Z\}$
 is a universal countable graph. It suffices to mention that for
 every pair $\bar x=(x_1,\dots, x_l)$, $\bar y=(y_1,\dots, y_p)$
  of finite sets from $\mathbb{Q}$ there exists a pattern whose white part
contains $\bar x$ and black part contains $\bar y$.
Finally, consider the closure $\bar Z$ of the set $Z$
in $\Bbb R$. We must prove that
the graph with $\Bbb R$
as the set of vertices and
$\{(x,y):x,y \in {\Bbb R}, |x-y|\in \bar Z\}$ as the set of edges is a universal graph.
Choose a pair $\bar x=(x_1,\dots, x_l)$, $\bar y=(y_1,\dots, y_p)$ of finite sets from
 $\mathbb{R}$
 and find a pattern $\gamma$ whose
 white part contains $\bar x$ and black part contains $\bar y$.
 The shift-invariance of the graph structure  follows from the construction.

2) In the case of a triangle-free graph we have only one additional
remark to our construction. As we have mentioned, the graph defined
in the induction base contains no
triangles; and, by the induction
hypothesis, no triangles appear when we define the sets $Z_i$, $i<n$.
Let us check that no triangles appear when we define the set $Z_n$.
Recall that we must consider only the white part of the pattern, because the
point $c$ is not joined by edges with the black part. But the slightly enlarged white part has no
edges by hypothesis, so the new edges do not produce triangles. As before, the
extension of the graph structure onto $\Bbb R$ is defined
by the closure $\bar Z$ of the set $Z$; since we have
chosen a sufficiently small open neighborhood of the point $c$,
the continuous graph inherits the absence of triangles.

 2. Now consider the case of a $K_s$-free universal countable graph
 for $s > 3$. The existence of a universal countable
 graph is a corollary of Fra{\"\i}ss\'e's axioms (one needs to check only the
 amalgamation axiom, see \cite{W}). But even in the
 countable case for $s >3$ there is no universal shift-invariant graph structure.
 More exactly, for a $K_s$-free universal graph
 there is no transitive action preserving graph-structure of the group $\Bbb Z$ on its vertices.
 Consider the case $s=4$. Assume that a shift-invariant universal $K_4$-free graph
 on the group $\Bbb Z$ does exist.
 Let $(0,a)$ be an edge. Choose $b$ such that $(b,0)$ is an edge, but $(b,a)$ and $(b,-a)$
 are not (this is possible by the universality). Then $(0,a+b)$ is not an edge
 (as well as $(b,-a)$), and hence the
 quadruple $(0,a,b,a+b)$ does not contain triangles (it is a quadrangle without diagonals).
 Hence, again by the universality, there exists $x$ joined with all points $0,a,b,a+b$.
 Then it is easy to check that the set $(0,a,x,x-b)$ is a 4-clique. We obtain a contradiction.
 In the case $s>4$ a contradiction may be obtained in a similar way, just start
 not from the edge $(0,a)$, but from some $(s-2)$-clique.

 The same claim is still true if we replace the group $\Bbb Z$ with an arbitrary abelian group.

 But the problem of constructing a universal continuous $K_s$-free graph for $s>3$
 without the requirement of shift-invariance is very easy.

 Let the set of vertices be $\Bbb R$. Again we define a pattern as a set of disjoint
 intervals with rational endpoints colored black and white. Let us enumerate all patterns
 as above. We will construct by induction a symmetric closed set $Z\subset
 \Bbb R \times \Bbb R$ with a nonempty interior, which will be the set of edges of our graph.
  As the induction base, we can choose $Z_1$ to be some square
  $[a,2a]\times [a,2a]$, $0<a$.
  At the $n$th step we consider the $n$th pattern $\gamma$ and fix the restriction of
  the set $\cup_{i=1}^{n-1}Z_i$ to the subgraph induced by the large segment $[-M_k,M_k]$,
 where $M_n=\max\{x:x\in \gamma_n \cup(-\gamma_n)\}+M_{n-1}+n+1$. Next we check
 whether there are cliques of size $n-1$ with white vertices of $\gamma_n$.
 If there are such cliques, we replace $n$ by $n+1$.
 If there are no such cliques, we add to the set
 $\cup_{i=1}^{n-1}Z_i$ the set $[M_n+1,M_n+2]\times \gamma_n^w$
 (where $\gamma_n^w$ stands for the white part of $\gamma_n$)
 and then symmetrize it in $\Bbb R \times \Bbb R$.
 It is easy to check that  after considering all patterns we get a
 topologically universal
 $K_s$-free graph.
 \end{proof}

We have proved that required topological and measurable universal
continuous  $K_s$-free graphs do exist, and this gives us examples of
${\frak S}^{\Bbb N}$-invariant measures on the space of adjacency
matrices.

The question which was posed in the first version of this paper
whether it is possible to construct a group-invariant structure
of a topologically universal graph for the compact group instead or $\mathbb R$
has an easy negative answer as the reviewer of the paper had mentioned.

Note that the concrete examples of  ${\frak S}^{\Bbb N}$-invariant measures on the set
of universal countable graphs that we have obtained here are new
and different from the Erd\"os--R\'enyi examples. For the construction
we use Theorems~3 and~4:

     Let $$dm(t)=\frac{1}{\sqrt {2\pi}} \exp\{-\frac{t^2}{2}\}dt$$
be the standard Gaussian measure on $\Bbb R$ and $\xi_1, \dots,
\xi_n, \dots$ be a sequence of independent random variables
each of which is distributed according to this Gaussian measure. Let
$E\equiv \{(t,s): |t-s|\in Z\}\subset {\Bbb R}^2$ where the set $Z$ was defined in the
proof of Theorem~4. Then the random $\{0;1\}$-matrices
           $$\{\chi_Z(\xi_i -\xi_j)\}_{i,j=1}^{\infty}$$
are, with probability 1, the adjacency matrices of universal (universal
triangle-free) graphs. In other words, the distribution of these
random matrices is a ${\frak S}^{\Bbb N}$-invariant measure
concentrated on the universal (triangle-free) graphs.
Of course, for the case of $K_s$-free graphs we also can choose the Gaussian
measure.
Instead of a Gaussian measure we can take any non-degenerate measure.
The choice of the set $Z$ ($E$) is not unique, as follows from the construction.

\section{Classification and the complete list of
invariant measures on the set of universal graphs}

As we have seen (Corollary~1), each measurable universal graph $(X,m,E)$ produces an invariant
measure on the set of universal countable graphs. Two questions arise:

1) When do two triples  $(X,m,E)$ and $(X',m',E')$ of universal measurable graphs produce the same
${\frak S}^{\Bbb N}$-invariant measures
$\mu_{\{X,m,E\}}$ and $\mu_{\{X',m',E'\}}$ on the set of universal graphs?

 Remark, that the list of invariant measures on the set of universal countable graphs
that are of type $\mu_{\{X,m,E\}}$
for some measurable universal graph $(X,m,E)$ (``randomization in vertices'') is not complete.
For example, it does not contain the Erd\'os--R\'enyi measure.
The second question is the following:

2)How to describe the complete list of ${\frak S}^{\Bbb N}$-invariant measures on the set of universal graphs?

We will give the answers to both questions.

\subsection{Classification of invariant measures obtained by randomization in vertices}

 The answer to the first question follows from a classification theorem of \cite{V},
which claims that two {\it pure} measurable symmetric functions of two variables
$f(x,y)$ and $f'(x',y')$ are isomorphic ($\Leftrightarrow$ there exists a
measure-preserving map $T:X\to X'$, $Tm=m'$, such that $f'(Tx,Ty)=f(x,y)$)
if and only if their matrix distributions coincide.
Recall that a pure symmetric function $f(\cdot,\cdot)$ of two variables is
a function for which
 the partition defined by the formula
 ($x\sim x_1 \Leftrightarrow f(x,y)=f(x_1,y)$ for almost all $y$)
 is the partition into separate points.

This property is true for a universal measurable graph. In our case,
the matrix distribution in the sense of $\cite{V}$ is just
the measure $\mu_{\{X,m,E\}}$.

\begin{theorem}
Two measurable universal graphs $(X,m,E)$ and $(X',m',E')$ produce
the same measure if and only if they are isomorphic in the following sense:
there exists a measure-preserving map $T:(X,m)\to (X'm')$ that sends
the set $E$ to $E'$:
$$(x,y)\in E\Leftrightarrow (Tx,Ty)\in E'.$$
\end{theorem}

Thus the measure $\mu_{\{X,m,E\}}$ is a complete isomorphism invariant
of measurable universal graphs. From this fact we immediately
obtain that our construction gives uncountably many different invariant
measures on the set of universal countable graphs, because even for a given
topologically universal graph $(X,E)$ we can vary the measure $m$ in
such a way that the measurable universal graphs $(X,m,E)$ are mutually
non-isomorphic for uncountably many measures $m$. It suffices to consider
$X=[0,1]$ with the Lebesgue measure $m$; then we can take uncountably many
symmetric sets $E\in X^2$
so that the measurable functions $x \mapsto m(E_x)$ for different choices of $E$
have nonequal distributions as measures on $[0,1]$; these
distributions are isomorphism invariants of the set $E$.

\subsection{Randomization in edges and description of the list of all invariant measures
on the universal homogeneous ($K_s$-free universal homogeneous) graphs.}

In order to obtain a description of the random countable universal graph ($K_s$-free-universal), or
in other words invariant ergodic
measure on the set of all countable  universal ($K_s$-free-universal) graphs),
we had considered the continuous universal graph and then chose by random the vertices
 of the countable graph. As we already have mentioned it is impossible to obtain
 the list of all possible invariant measures on the set of universal graphs with this procedure --- this
 ``randomization in vertices'' only,  but another source of the randomness is "randomness in the edges",
 which allow to obtain whole list of invariant measures on the set of universal countable graphs.
 Below we explain what does this mean. But for $K_s$-free universal graphs we
 do not need the randomization
 in edges. We see below that it is enough to use randomization in vertices only.

Instead of the measurable graph $(X,m,E)$ we consider more general object.

\begin{definition}
Generalized measurable continuous graph is the triple $(X,m,\omega)$, where $(X,m)$ is the standard measure space with
continuous normalized measure (Lebesgue space) and $\omega$ is any symmetric measurable function on the space
$(X\times X, m\times m)$ with values $\omega(x,y)\in[0,1]$
\end{definition}

\begin{remark}
In the case the function $\omega$ takes value $\{0;1\}$ - (we will call this case -deterministic in edges)
the subset $E=\{(x,y):\omega(x,y)=1\}$ gives
the definition of the measurable graph $(X,m,E)$ in the sense of paragraph 2
\end{remark}

The interpretation of the value of the function $\omega$ at the
point $(x,y)$ is a probability that $(x,y)$ is the edge of the
continuous generalized graph $(X,m,\omega)$.

\subsubsection{The list of all invariant measures for the case of universal homogeneous graphs}

\textbf{Construction}.
{\it This is two-step randomization.
 Suppose we have continuous generalized graph $(X,m,\omega)$; it produces a
measure on the space of all countable graphs (or produces a random graph) as follows:

{\rm1.} we choose the set of vertices $\{x_k\}_{k=1}^{\infty}$ as a sequence of independent random points from $X$
with respect to measure $m$ ("randomization in vertices");

{\rm2.} for each chosen pairs of the vertices $(x_i,x_j)$ we define whether this
pair is edge or not in our random graph
independently (over all i,j) with probabilities $\omega(x_i,x_j), 1-\omega(x_i,x_j)$.}

For any generalized graph $(X,m,\omega)$ this construction gives the measure, on the space of
countable graphs
which we denote as $M(X,m,\omega)$; will say that this measure is generated by
generalized continuous graph $(X,m,\omega)$.
So, we have the map from the set of the generalized graphs to the set of
measures on the space of countable graphs
(or its adjacent matrices).

This double randomization can be considered as randomization of the probabilities --- first step,
and consequent realization of those probabilities ---
the second step. Such tool is typical for the theory of the random walk in random environment.

Let us give the precise formula for the measure  $M(X,m,\omega)$ of the cylindric sets. Suppose $A=\{a_{i,j}\}, i,j=1 \dots n$
is $(0-1)$-matrix of order $n$, and $C_A$ is a cylindric set of all infinite $(0-1)$-matrices, which has the matrix $A$ as submatrix
on the NW-corner. Then the value of the measure $M(X,m,\omega)\equiv M$ on the cylinder $C_A$ is
   $$M(C_A)=\int_{X^n} \prod_{(i,j): a_{i,j}=1} \omega(x_i,x_j) \prod_{(k,r): a_{k,r}=0}(1-\omega(x_k,x_r))dm(x_1)\dots dm(x_n).$$
This is a generalization of the formula given after theorem 3
in deterministic in edges case.

The following fact is evident from the definition
\begin{theorem}
The measure $M(X,m,\omega)$ is ${\frak S}^{\Bbb N}$-invariant and ergodic.
\end{theorem}

Now we formulate in the convenient for the case of graphs
form of the theorem by D.~Aldous (\cite{A}) which gives the description
of  ${\frak S}^{\Bbb N}$-invariant measures on the space infinite matrices.

\begin{theorem}
Each ergodic ${\frak S}^{\Bbb N}$-invariant measures of the
set of symmetric with zero diagonal infinite $\{0;1\}$-matrices
is generated by generalized graph in the framework of the construction above.
\end{theorem}

Consequently, each ergodic ${\frak S}^{\Bbb N}$-invariant measures on the set of all countable graphs
can be obtained by or construction above. We will not discuss here the proof of Aldous theorem but remark that
the second author will present elsewhere an alternative approach ("ergodic method") to
this theorem and will give a new proof of it
(see also discussion in \cite{V}).

Remark also, that randomization in edges (when exist) cannot be reduced to randomization in vertices:
more exactly, the resulting measure $M(X,m,\omega)$ can not be obtained as a measure corresponding to
the measurable graph $(X,m,E)$ \footnote{It corresponds in a sense to the generalized
 function $\omega$ (see \cite{V})}.

Now we must formulate the special condition on the generalized graph when the construction above gives
the measure {\it on the universal ($K_s$-free universal) homogeneous graph}. The condition
is similar to the condition of the theorem 2.
We will formulate it only for universal graphs: for $K_s$-free universal graphs $s>2$ we
do not need such notion (see below).

\begin{definition}
The generalized measurable graph  $(X,m,\omega)$ is called universal if for all
natural $n,m$ and almost all pairs of
the sets $(x_1, \dots, x_n), (y_1, \dots, y_m)$ from $X$ the
following is true
$$m\{z \in X:\prod_{i,j}\omega(x_i,z)(1-\omega(y_j,z))>0\}>0.$$
\end{definition}

Here were restrict ourselves with the following important remark. The first step
of the construction gives us
a family  Bernoulli measures on the {0-1}-matrices (all entries are independent
but have in general different distributions)
 and the resulting measure on the space of $0-1$-matrices is the average
 (more exactly barycenter) of the those Bernoulli measures.

The explicit construction of the universal ($K_s$-free universal) generalized graph can be done
(even simpler) as for
universal topological graph in the section 3 which is of course partial case. We will not consider this,
but emphasize that the constructions of the section 3 have more instructive character than construction of
generalized universal graph.

Now we can give a list of all ${\frak S}^{\Bbb N}$-invariant ergodic measures on the space
of all universal countable graphs
(or on the set of corresponding ${0;1}$ matrices), which is the mail goal of this section.
This is the corollary of the previous theorems:
\begin{theorem}
Each ergodic ${\frak S}^{\Bbb N}$-invariant measures on the set of all
countable universal homogeneous graphs is generated by
the construction above with generalized measurable universal graph.
\end{theorem}

\subsubsection{The list of measures for $K_s$-free universal homogeneous graphs.}

\begin{theorem}
The list of all invariant ergodic measures on the $K_s$-free (for $s>2$) countable homogeneous graphs
is given by randomization in vertices only, e.f. by the construction of theorem 3 of the subsection 2.2.
In other words, in order to obtain whole list it is enough to use only the randomization in vertices
of measurable graph and no randomization in edges.
\end{theorem}
We outline the simple proof of this fact. Let for simplicity consider the case $s=3$; the general case is
analogous. Assume the contrary, and suppose that there exist
a set of positive measure $F \subset X^2$ such that
 for $(a,b)\in F, 0<\epsilon_1 <\omega (a,b)<\epsilon_2 <1$. Then on the one hand, $a$ and $b$ have almost
 no common neighbor for almost all $(a,b)\in B$ (since the graph must be triangle-free), on the other hand
 $a$ and $b$ for almost all $(a,b)\in B$ must have a common neighbor with positive probability
 (since the graph must be universal, and if we do not take the edge a-b, then these vertices should
 have some common neighbor).

\section{Some problems and comments}

\smallskip\noindent 1. \textbf{The distribution of the entries of random adjacency}. A very important
question is to characterize, for an arbitrary invariant measure
$M(X,m,\omega)$ (see the previous section), the finite-dimensional
distribution of the entries with respect to this measure. Because of the ${\frak S}^{\Bbb N}$-invariance
of the measure, the finite-dimensional distributions are
${\frak S}_n$-invariant; consequently, they are concentrated on bunches of
orbits of these groups and decompose into positive combinations of orbits
of the group ${\frak S}_n$ in the space of matrices $M_n(0,1)$. Perhaps,
because of the ergodicity of the measure $M(X,m,\omega)$,
the finite-dimensional distributions must concentrate near
one or several typical orbits at short distances from one another.
This is an analog of the Law of Large Numbers.
How to characterize these orbits? The answer could be useful for the solution of the problem by Cherlin
(see introduction). The structure and the asymptotic size of these
orbits is an interesting characteristic of universal graphs and the measures.

\smallskip\noindent
2. \textbf{Uniqueness of
measurable universal graphs.}
When discussing the definitions of universality above, the following question
naturally arises: under what conditions is the set $E$ which defines
the graph structure on the standard Borel (or standard measure) space
of universal Borel or measurable graphs unique up to
isomorphism? In the case of countable universal graphs, the
``back and forth'' method allows one to prove the uniqueness of the universal graph.
Equivalently, the question above is as follows: when does
the Borel or measurable version of the ``back and forth''
method work? The same question can be solved positively for metric spaces:
as proved by Urysohn, there exists a unique (up to isomorphism)
universal Polish space.
It is interesting to have link with model theory in which one consider finite or countable situation.
When uniqueness takes place in the continuous case?

\smallskip\noindent 3. \textbf{Approximation.} In our construction we obtained a
continuous graph
with $\Bbb R$ as the set of vertices as the completion of a graph with
the set of vertices $\Bbb Q$.  Of course, in that case we could define the graph directly,
avoiding approximation.
But it is interesting whether in the general situation of model theory it is possible to
consider a ``completion'' of countable models. More exactly, how to formulate
Fra{\"\i}ss\'e's axioms
for the Borel or measurable case (with separability conditions)
in order to obtain it as the projective or another limit of  the finite theory?
A very good example of a positive solution of such a problem is, of course, the theory
of universal metric spaces.

\smallskip\noindent 4. \textbf{Link to the Urysohn space.} In this
sense, the Urysohn space is of special interest. We will consider it from this
point of view elsewhere. Here we mention only that the  Urysohn space
$\Bbb U$ plays the role of a ``Borel universal object'' (or
topologically universal object) for the rational or integer universal
metric space. Any Borel probability measure $m$ on this space
defines a ${\frak S}^{\Bbb N}$-invariant measure $\mu$ on the
space $\cal R$ of {\it distance real matrices} which are
universal\footnote{The universality of a distance matrix
means that the completion of $\Bbb N$
with respect to the corresponding metric is isometric to the Urysohn space,
see \cite{V1}} with probability one. The similarity between the theory of the Urysohn space and the
example of Section~3 above can be illustrated by the result of
\cite{CV} where the Urysohn space was realized as the completion of the
real line with respect to a universal shift-invariant metric.

\end{document}